\documentclass[aos,MSNbibl,nameyear,dvips]{arximspdf}
\usepackage{subenv}

\doi{10.1214/13-AOS1083} 
\volume{41}
\issue{1}
\pubyear{2013}
\firstpage{323}
\lastpage{341}

\makeatletter
\newcommand{\rrVert}{\Vert}
\newcommand{\rrvert}{\vert}
\newcommand{\llVert}{\Vert}
\newcommand{\llvert}{\vert}

\newcommand{\cal}{\mathcal}

\newproclaim{Def}{Definition}
\newtheorem{theorem}{Theorem}

\newcommand{\ceps}{{\cal E}}
\newcommand{\ca}{{\cal A}}

\newcommand{\pr}{\mathrm{P}}

\newcommand{\st}{\mid}
\newcommand{\E}{E}
\newcommand{\into}{\rightarrow}
\newcommand{\NatN}{\mathbb N}

\newcommand{\CompN}{\mathbb C}
\newcommand{\QuartN}{\mathbb H}
\newcommand{\OctN}{\mathbb O}
\newcommand{\RealN}{\mathbb R}

\makeatother

\begin{document}
\begin{frontmatter}

\title{Fiducial theory and optimal inference}
\runtitle{Fiducial theory}

\begin{aug}
\author[A]{\fnms{Gunnar} \snm{Taraldsen}\corref{}\ead[label=e1]{Gunnar.Taraldsen@sintef.no}}
\and
\author[B]{\fnms{Bo Henry} \snm{Lindqvist}\ead[label=e2]{bo@math.ntnu.no}}
\runauthor{G. Taraldsen and B. H. Lindqvist}
\affiliation{SINTEF Information and Communication Technology and
Norwegian~University~of Science and Technology}
\address[A]{SINTEF ICT\\
P.O. Box 4760 Sluppen\\
NO-7465 Trondheim\\
Norway\\
\printead{e1}}
\address[B]{Department of Mathematical Sciences\\
Norwegian University of Science\\
\quad and Technology\\
NO-7491 Trondheim\\
Norway\\
\printead{e2}} 
\end{aug}

\received{\smonth{6} \syear{2012}}
\revised{\smonth{1} \syear{2013}}

%
\begin{abstract}
It is shown that the fiducial distribution in a group model, or more
generally a quasigroup model, determines the optimal equivariant
frequentist inference procedures. The proof does not rely on existence
of invariant measures, and generalizes results corresponding to the
choice of the right Haar measure as a Bayesian prior.
Classical and more recent examples show that fiducial arguments can be
used to give good candidates for exact or approximate confidence
distributions. It is here suggested that the fiducial algorithm can be
considered as an alternative to the Bayesian algorithm for the
construction of good
frequentist inference procedures more generally.
\end{abstract}

%
\begin{keyword}[class=AMS]
\kwd[Primary ]{62C05}
\kwd[; secondary ]{62A01}
\kwd{62F10}
\kwd{62F25}
\kwd{20N05}
\end{keyword}
\begin{keyword}
\kwd{Fiducial}
\kwd{group model}
\kwd{invariance}
\kwd{risk}
\kwd{frequentist inference}
\end{keyword}

\end{frontmatter}

\section{Introduction}\label{intro}

Fiducial theory was introduced by \citet{Fisher30} to avoid the
problems related to the choice of a prior distribution.
Fiducial inference has not gained much popularity as such,
but the related theory has been historically influential
[\citet{Efron98}],
and is still important in the current flow of statistical developments
[\citet{Efron06fiducial,LidongHannigIyer08OneWay,GhoshReidFraser10CondInference,FraserReidMarrasYi10fiducial,WangHannigIyer12gumfiducial}].
Fisher's own view on fiducial inference evolved over the years as can
be inferred from a reading of his initial
[Fisher (\citeyear{Fisher30,Fisher35Fiducial})] and more final formulation
of the theory [\citet{FISHER}].
He was in particular very positive to the developments by
Fraser (\citeyear{Fraser61Fiducial,Fraser61FiducialInvariance,Fraser62fiducial,Fraser64fiducial}),
and we most certainly share this point of view.
Fraser (\citeyear{FRASER,Fraser79inferenceAndLinear}) develops the theory and uses
the label \textit{structural inference} for this.
A strongly related theory was presented under the label of
\textit{functional models} by \citet{Bunke75functional}
and \citet{DawidStone82}.
The term \textit{fiducial} will here be used more generally so that
it includes structural, functional, and the original fiducial
arguments given by Fisher.

The original idea of Fisher was to obtain the fiducial distribution
directly from the cumulative distribution,
but this line of argument runs into problems when similar
arguments are tried in the multivariate case.
The definition used here is based on the solution of a fiducial equation,
and is in this sense similar to the approaches of
\citet{FRASER}, \citet{DawidStone82} and
Hannig (\citeyear{Hannig09fiducial,Hannig12fiducial}).
A~more precise definition of the term \textit{fiducial model}
as used here is given in
Section~\ref{sOptimal} in Definition~\ref{defFidMod}.
A brief review of alternative, but strongly related definitions
found in the literature is given in the final Section~\ref{sDiscussion}.


Let $l = \gamma(\theta,a)$ denote the realized loss for an action $a
\in\Omega_A$ given the model parameter $\theta\in\Omega_\Theta$. Let
$\Omega_X$ be the sample space equipped with the $\sigma$-field
${\cal
E}_X$ of events. The risk $\rho$ of a decision rule $\delta\dvtx  \Omega_X
\into\Omega_A$ is by definition equal to the expected value $\rho=
\E^\theta\gamma(\theta, \delta(X))$. This is determined by the
statistical model given by the family $\{(\Omega_X, {\cal E}_X,
\pr_X^\theta) \st\theta\in\Omega_\Theta\}$ of probability spaces.

Consider now the more special case where
$\Omega_X = \Omega_\Theta= \Omega_A = G$,
possibly after a suitable change of variables.
Assume that $G$ is
a measurable quasigroup with a unit $e$,
and product $(g_1, g_2) \mapsto g_1 g_2$ written like ordinary
multiplication [\citet{Smith06loop}].
This includes the more common case of a group,
but it is more general since the associative law is not assumed to hold.
Assume furthermore that
$X \sim\theta U$ conditionally on $\Theta= \theta$ and that
the law of $(U \st\Theta= \theta)$
does not depend on $\theta$.

This gives an example of
a fiducial model for the statistical model
as defined more generally
on page \pageref{defFidMod}.
The fiducial distribution is
obtained by solving the fiducial equation $x = \theta u$
for $\theta$ when $u$ is sampled from $\pr_U^\theta$.
Existence and uniqueness is ensured since
$G$ is a quasi-group.
A variable $\Theta^x$
is uniquely determined by $x = \Theta^x U$.
The fiducial distribution is then the conditional law of $\Theta^x$ given
$\Theta= \theta$.

Assume that the loss is invariant
in the sense that
$\gamma(\theta, a) = \gamma(g \theta, g a)$,
and that the decision rule is equivariant in the sense
that $\delta(g x) = g \delta(x)$.
The assumptions ensure the validity of
the following calculation:
%
\begin{subequation}
\label{eqRiskCalc}
%
\begin{eqnarray}
\label{eqRisk1}
\rho& = & \E^\theta\gamma\bigl(\theta, \delta(X)\bigr) \hspace*{22.86pt}\qquad\mathrm
{Risk}
\\[-2pt]
& = & \E^\theta\gamma\bigl(\theta, \delta(\theta U)\bigr) \hspace*{16.74pt}\qquad
\mbox{Fiducial model for } \pr_X^\theta
\\[-2pt]
& = & \E^\theta\gamma\bigl(\theta, \theta\delta(U)\bigr) \hspace*{16.95pt}\qquad\mbox
{Equivariance of $\delta$}
\\[-2pt]
\label{eqRisk4}
& = & \E^\theta\gamma\bigl(e, \delta(U)\bigr) \hspace*{23pt}\qquad\mbox{Invariance of
$\gamma$}
\\[-2pt]
& = & \E^\theta\gamma\bigl(\Theta^x, \Theta^x
\delta(U)\bigr) \qquad \mbox{Invariance of $\gamma$}
\\[-2pt]
& = & \E^\theta\gamma\bigl(\Theta^x, \delta\bigl(
\Theta^x U\bigr)\bigr) \hspace*{1pt}\qquad \mbox{Equivariance of $\delta$}
\\[-2pt]
\label{eqQFid}
& = & \E^\theta\gamma\bigl(\Theta^x, \delta(x)\bigr) \hspace*{17.74pt}\qquad
\mbox{Fiducial equation}.
\end{eqnarray}
\end{subequation}
%
%
A variation of the above argument gives
that $\Theta^x$ can be replaced by
$x U_r^{-1}$ in the conclusion.
In the group case the law of
$\Theta^x$ will coincide with the law of\vadjust{\goodbreak}
$x U_r^{-1}$,
but in general not since the defining equation $e=U U_r^{-1}$
of the right inverse $U_r^{-1}$ does not provide the solution of the
fiducial equation.
It follows from this that an optimal equivariant rule,
if it exists,
is determined by the fiducial distribution of $\Theta^x$
or by the distribution of $x U_r^{-1}$ from the right inverse.
The first part of the claims in the abstract has hence been established.

It is hoped that the reader can appreciate the simplicity and
consequence of the calculation given in equation (\ref{eqRiskCalc}),
but it could also be considered to be essentially Greek.
The required theory of decisions and fiducial theory will be
explained in some more detail in Section~\ref{sOptimal},
and examples are presented in Section~\ref{sExamples}.
The presentation is essentially as given in standard textbooks
[\citet
{LehmannCasella98Estimation,LehmannRomano05testing,SCHERVISH,BERGER,StuartOrdArnold99kendal2a}],
but with the simplifications given by a fiducial model.
The monographs by \citet{Eaton89groupstat} and
\citet{Wijsman90statgroups} are recommended as excellent sources
for theory and examples beyond the standard textbooks. 
%
%

The presentation in the following will be mostly restricted to the
group case,
but it will be more general than the previous in the sense that
the assumption of equality of the involved spaces will be abandoned.
It will be more general than standard theory since, as above,
the arguments will not depend on existence of invariant measures.

\section{Optimal inference}
\label{sOptimal}

Consider the case where the loss of an action $a \in\Omega_A$
is of the form $l = \gamma(\theta, a)$
corresponding to a statistical model $\{\pr_X^\theta\st\theta\in
\Omega_\Theta\}$.
It is here assumed that the model parameter $\Theta$ is a $\sigma$-finite
random quantity and this and all other
random quantities are defined based
on the underlying conditional probability space $(\Omega, \ceps, \pr)$
as explained by \citet{TaraldsenLindqvist10ImproperPriors}.
This means in particular that
$\pr_X^\theta(B) = \pr(X \in B \st\Theta= \theta)$,
and $X\dvtx  \Omega\into\Omega_X$, $\Theta\dvtx  \Omega\into\Omega_\Theta$
are measurable functions.
It means also that all expectations that occur are defined by
integration over $\Omega$.
As an example
$\E(\phi(Z) \st T=t) = \int\phi(Z(\omega)) \pr^t (d\omega)$
by definition.
It is here assumed that
$\phi\dvtx  \Omega_Z \into\RealN$,
$Z\dvtx  \Omega\into\Omega_Z$, and
$T\dvtx  \Omega\into\Omega_T$ are measurable.
The conditional law $\pr^t$ is well defined if $\pr_T$ is
$\sigma$-finite.
The consequence
$\E(\phi(Z) \st T=t) = \int\phi(z) \pr_Z^t (dz)$
is a theorem.

The law $\pr_\Theta$ of $\Theta$ is not assumed known and is not
needed in the arguments which follow.
The reason for the assumption of existence of $\Theta$, $X$,
and indeed any random quantity involved in the arguments,
as functions defined on the conditional probability space
$(\Omega, \pr, \ceps)$
is as in the formulation of probability theory given by
\citet{KOLMOGOROV}:
any collection of random quantities gives a new random quantity with a
well-defined law,
and measurable functions of random quantities give new random quantities.
The theory is completely based on the underlying abstract space $\Omega$.

A \textit{group invariant problem} is given by a group $G$ that
has a transformation group action on the sample space $\Omega_X$,
the model parameter space $\Omega_\Theta$,
and the action space $\Omega_A$.\vadjust{\goodbreak}
The problem is group invariant if
$\pr_{gX}^\theta= \pr_X^{g \theta}$
and $\gamma(g \theta,g a) = \gamma(\theta, a)$.
An inference rule $\delta$ with a corresponding
action $A=\delta(X)$ is \textit{equivariant} if
$\delta(g x) = g \delta(x)$.
The restriction to the class of equivariant actions can be
interpreted as a consistency requirement:
an observation $x$ from $\pr_X^\theta$
corresponds to an observation
$g x$ from $\pr_X^{g \theta}$.
The two corresponding problems are formally identical and the
use of an equivariant action ensures consistency.

The problem considered here is to
determine an equivariant $\delta$ such that the risk
%
\begin{equation}
\label{eqRisk} \rho= \E^\theta\gamma\bigl(\Theta, \delta(X)\bigr) = \E
\bigl(\gamma\bigl(\Theta, \delta(X)\bigr) \st\Theta= \theta\bigr)
\end{equation}
is minimized.
It will be assumed that $G = \Omega_\Theta$ with the action
given by the group multiplication $g \theta$ directly.
The \textit{orbit} of $x$ in $\Omega_X$ is defined by
$Gx = \{g x \st g \in G\}$,
and likewise for orbits in $\Omega_\Theta$ and $\Omega_A$.
The action of $G$ is \textit{free} on $\Omega_X$ if the mapping
$g \mapsto g x$ is injective for all $x$.
The group action is \textit{transitive} on $\Omega_X$ if
$G x = \Omega_X$.
If the group action is both transitive and free,
then it is said to be \textit{regular} and the corresponding space
is then a principal homogeneous space for $G$.
It follows in particular that the model parameter space
$\Omega_\Theta$ is a principal homogeneous space for $G$,
but there has also been an identification of the identity element $e$
in $\Omega_\Theta$.

Let $U$ be a random quantity such that
\mbox{$\pr(U \in A \st\Theta= \theta) = \pr(X \in A \st\Theta=e)$}
holds identically for all $A$ and $\theta$.
It follows then that
%
\begin{equation}
\label{eqGFid2} (X \st\Theta= \theta) \sim(\theta U \st\Theta= \theta)
\end{equation}
since the group invariance of the statistical model justifies
$\pr_X^\theta= \pr_{\theta X}^e = \pr_{\theta U}^\theta$.
This construction proves that
$(U, \chi)$ with
%
\begin{equation}
\label{eqGroupMod} \chi(u,\theta) = \theta u
\end{equation}
is a fiducial model for $\pr_X^\theta$.
The concept of a fiducial model is defined as follows.
%
\begin{Def}[(Fiducial model)]
\label{defFidMod}
Let $\Theta$ be a $\sigma$-finite random quantity.
A fiducial model $(U, \zeta)$ is given
by a random quantity $U$ and a measurable function
$\zeta\dvtx  \Omega_U \times\Omega_\Theta\into\Omega_Z$.
This is a fiducial model for the
statistical model $\{\pr_Z^\theta\st\theta\in\Omega_\Theta\}$ if
%
\begin{equation}
\bigl(\zeta(U,\Theta) \st\Theta= \theta\bigr) \sim(Z \st\Theta=
\theta).
\end{equation}
\end{Def}
The notation
$(W_1 \st\Theta= \theta) \sim(W_2 \st\Theta= \theta)$
means that
$\pr_{W_1}^\theta= \pr_{W_2}^\theta$ so Definition~\ref{defFidMod}
can be compared with equation (\ref{eqGFid2}).
It is allowed in the above that
$\pr_U^\theta$ does depend on $\theta$.
Interesting examples where this occurs are discussed by
\citet{Fraser79inferenceAndLinear} in the form of dependence on
shape parameters in addition to pure group parameters.
In the following it will, however, be assumed throughout that
the fiducial model is \textit{conventional} in the sense that
$\pr_U^\theta$ does not depend on $\theta$.\vadjust{\goodbreak}

It is important to notice that many different fiducial models are possible
for a given statistical model.
A fiducial model provides a different
basis for statistical inference than a statistical model.
The choice of a particular fiducial model can be compared
with the choice of a Bayesian prior together with a statistical model.
Fiducial inference is then initially different from frequentist and
Bayesian inference since the inferential basis is
given by a fiducial model which is assumed known.
Fiducial inference as such will not be considered here,
but the corresponding fiducial algorithms will be used as vehicles for
the construction
of frequentist procedures.

A fiducial model $(U, \zeta)$ is \textit{simple} if
the fiducial equation $\zeta(u,\theta) = z$ has a unique
solution $\theta^z (u)$ when solved for $\theta$ for all $u,z$.
In the simple and conventional case the fiducial distribution is defined
as the distribution of $\Theta^z = \theta^z (U)$ conditional on
$\Theta= \theta$.
%
\begin{Def}[(Fiducial distribution in the simple and conventional case)]
\label{defFid1}
Let $(U, \zeta)$ be a conventional simple fiducial model.
Define the random quantity $\Theta^z$ by
$z = \zeta(U, \Theta^z)$.
The fiducial distribution is the conditional law of
$\Theta^z$ given $\Theta= \theta$.
\end{Def}

The fiducial model $(U, \chi)$ given by equation (\ref{eqGroupMod})
is simple if and only if
$\Omega_X$ is a principal homogeneous space for $G$.
In this case, by the choice of a unit element in $\Omega_X$,
the identification $G=\Omega_\Theta=\Omega_X$ can be done.
It follows then that $\theta^x (u) = x u^{-1}$ is the unique solution
of $x = \theta u$,
and the fiducial distribution is the conditional distribution
of $\Theta^x = x U^{-1}$ as it appears in equation (\ref{eqQFid}).


The remainder of this section will be on the analysis
of the group model by means of the constructed fiducial model given by
equation (\ref{eqGroupMod}) in the
case where $\Omega_X$ is not assumed to be a principal homogeneous
space for $G$.
The aim is to determine an equivariant
inference rule $\delta$ so that the risk given by
equation (\ref{eqRisk}) is minimized.
A definition of a fiducial distribution will also
be presented for this group case.
The resulting distribution coincides with the distribution described
with many more examples,
explicit calculation of densities,
and illustrative figures by
Fraser (\citeyear{FRASER,Fraser79inferenceAndLinear}).

A first observation is that the calculations given by
equations (\ref{eqRisk1})--(\ref{eqRisk4}) are valid and the risk is
given by
$\rho= \E(\gamma(e,\delta(U)) \st\Theta= \theta)$.
The construction of the fiducial model has hence given a simple proof
that gives
that the risk does not depend on the model parameter since
$\pr_U^\theta$ does not depend on $\theta$.

Let $Y = \phi(X)$ be an invariant statistic in the sense
that $\phi(\theta x) = \phi(x)$ for all $\theta,x$.
This is equivalent with the requirement that $\phi$ is constant on
all orbits in the sample space $\Omega_X$.
The fiducial model in equation (\ref{eqGroupMod}) gives
that $Y = \phi(X) \sim\phi(\Theta U) = \phi(U)$
conditionally on $\Theta= \theta$.
The conclusion is that $\pr_Y^\theta$ does not depend on
$\theta$ and has a known distribution.
This proves that an invariant statistic $Y$ is ancillary.\vadjust{\goodbreak}

Assume furthermore that $Y = \phi(X)$ is a \textit{maximal invariant statistic}.
This means that the family of level sets of $\phi$ coincides with
the family of orbits in $\Omega_X$.
Let $x$ be given and assume that
$y = \phi(x) = \phi(u)$.
The maximality ensures that $x \in G u = \Omega_\Theta u$,
so $x = \theta^x u$ for some $\theta^x$.
This $\theta^x$ will be unique if $G$ acts freely on $\Omega_X$,
but here it will more generally be assumed that
$\theta^x$ is determined by the choice of a measurable selection.
The measurable selection theorem [\citet{CASTAING}]
ensures existence under mild conditions.
The fiducial distribution of the corresponding variable
$\Theta^x$ can be described as follows.
%
\begin{Def}[(Fiducial distribution in the group case)]
\label{defFid2}
Let
$u$ be a sample from the distribution of
$(U \st\Theta=\theta, \phi(U) = \phi(x))$ where $\phi$ is a
maximal invariant.
Let $\theta^x$ be a measurable selection solution
of $x = \theta^x u$.
This $\theta^x$ is a sample from a fiducial distribution.
\end{Def}
The solution $\theta^x$ exists since $y = \phi(x) = \phi(u)$ ensures
that $x$ and $u$ are on the same orbit.
Definition~\ref{defFid2} is a special case of
Definition~\ref{defFid1} if $\Omega_X$ is a principal homogeneous
space for $G$,
and the definitions are hence consistent.
It is possible to define a fiducial distribution for more general cases.
One version is as presented by \citet{Hannig09fiducial},
but there are also other possibilities available.
This will not be discussed further here since the given definitions
of the fiducial distribution are sufficient for the purposes
in this paper.

Let $Y = \phi(X)$ be a maximal invariant statistic.
The calculation that gave equation (\ref{eqRisk4}) can
now be continued to give
%
\begin{equation}
\rho = \int \bigl[\E^{\theta,y} \gamma\bigl(e,\delta(U)\bigr) \bigr]
\pr^\theta_Y (dy). 
\end{equation}
The expression $[\cdot]$ does only depend on $y$.
The optimal rule $\delta$, if it exists,
is found by minimization for each given $y$.
Assume that $x$ is such that $y=\phi(x)$.
It follows then that
%
\begin{equation}
\E^{\theta,y} \gamma\bigl(e,\delta(U)\bigr) = \E^{\theta,y} \gamma\bigl(
\Theta^x, \Theta^x \delta(U)\bigr) = \E^{\theta,y}
\gamma\bigl(\Theta^x, \delta(x)\bigr)
\end{equation}
and the optimal rule $\delta$ is determined by
the fiducial distribution of $\Theta^x$.
The variable $\Theta^x$ is defined as
a measurable selection solution of $x = \Theta^x U$.
This result can be summarized as the main technical result in this paper.
%
\begin{theorem}
\label{theo1}
The risk of an equivariant rule in a group invariant problem is
determined by
a fiducial distribution if the model parameter space is
a principal homogeneous space for the group.
\end{theorem}
It should be noted that the statement assumes existence of a fiducial
distribution as described above,
but uniqueness of a fiducial distribution is not assumed.
Optimal inference procedures
are determined by the fiducial distribution regardless
of the choice of a measurable selection for the determination of a
fiducial distribution.\vadjust{\goodbreak}
The optimal $\delta$ is found, if it exists, as the minimizer $\delta
(x)$ of the expression
%
\begin{equation}
\E^{\theta,y} \gamma\bigl(\Theta^x, \delta(x)\bigr),
\end{equation}
where the conditional distribution of $\Theta^x$ is a fiducial distribution.

Theorem~\ref{theo1} generalizes directly to the larger class of
randomized equivariant actions.
This is obtained by a replacement of the equivariant action
$\delta(X)$ by the randomized equivariant action
$\delta(X, V) = \delta(\theta U, V)$ in the calculations.
It is here assumed that $U$ and $V$ are conditionally
independent in the sense that
$\pr_{U,V}^\theta(du,dv) = \pr_U^\theta(du) \pr_V^\theta(dv)$,
and both conditional distributions do not depend on $\theta$.
The equivariance is defined by the identity
$\delta(g X, V) = g \delta(X, V)$.

A randomized action corresponds to the assignment of a probability
measure on the action space $\Omega_A$. The set of randomized actions
is hence always a convex set, and this gives theoretical advantages to
the problem of minimization of the risk. If, however, the loss function
$l (\theta, a)$ is convex on $\Omega_A$ for each $\theta$, then the
Jensen inequality gives that it is sufficient to consider nonrandom
actions [\citet{LehmannCasella98Estimation}, page 48].

Theorem~\ref{theo1} generalizes also directly to the case where
$G$ is only assumed to act transitively on $\Omega_\Theta$.
The construction is as above,
and starts with fixing a $\theta_0$ and the construction of
a random variable $U$ such that $\pr_U^\theta= \pr_X^{\theta_0}$.
All the arguments given above can then be repeated with
$G$ playing the role of a new and possibly larger parameter space.
The result is then first a fiducial distribution on $G$,
but this is pushed forward to a fiducial distribution on
$\Omega_\Theta$ by the mapping $g \mapsto g \theta_0$.


It is known that the fiducial coincides with the posterior from a right
Haar prior,
and for these cases Theorem~\ref{theo1} is a known result with
the posterior used in the formulation instead.
There are, however, groups where no Haar prior exists, and in this
case Theorem~\ref{theo1} and its extensions given by the above
comments is a novel result.
The derivation given in the \hyperref[intro]{Introduction} also gives a similar result
in the more general case of a quasi-group,
and the existence of invariant measures is then also not automatic.

\section{Examples}
\label{sExamples}

The examples presented next are
chosen to illustrate the concepts.
Many more examples and thorough discussions are found in
the previously quoted textbooks and monographs.
A complete treatment of the given examples---including in particular
simulation studies of the resulting procedures---will not be pursued since
this would tend to take attention away from the main issue.
The purpose is simply to indicate the usefulness of fiducial theory.

\subsection{The Bernoulli distribution}
\label{ssBernoulli}

A random sample of size n from the Bernoulli distribution provides an
example where the results related to
Theorem~\ref{theo1} cannot be applied directly.
Fiducial theory can, nonetheless, be used to obtain optimal inference.\vadjust{\goodbreak}

The largest possible group $G$ equals $\{e, g_1\}$ corresponding
to the group of permutations of the set $\{0,1\}$.
The action on $\Omega_\Theta= (0,1)$ is determined by $g_1 p = 1-p$,
and the set of orbits in the parameter space is uncountable.
The conclusion is that conditioning on the maximal invariant as
in the arguments leading to Theorem~\ref{theo1} does not
provide any essential simplification of the problem.

This example is, however, very important from the point of view that
fiducial distributions can still provide optimal procedures.
\citet{Blank56binomialoptimal} has constructed a randomized
most powerful unbiased confidence interval,
and this is related to
a fiducial distribution [\citet
{Anscombe48fiducial}, Stevens (\citeyear{Stevens50binomial,Stevens57binomial})].

The empirical mean is the unbiased estimator of $p$ with minimum variance.
It can, however, be argued that neither unbiasedness nor minimum
variance are natural concepts in this particular case.
The parameter space $\Omega_\Theta$ can alternatively
be identified with the circular arc
$\{(\sqrt{p},\sqrt{q}) \st p,q>0, p+q=1\}$ in the
$(\sqrt{p},\sqrt{q})$-plane.
This has the advantage that the Fisher information
metric distance between two
distributions in this parametric family equals the
distance along the arc [\citet
{Rao45,AtkinsonMitchell81statdistance,Amari90diffgeostat}].
The distance squared provides a loss that is invariant
with respect to $G$.
A natural task is to investigate on existence
of an optimal equivariant estimator
of $p$ with respect to the distance squared on the arc.
A reasonable candidate arises from the previously referenced fiducial
distribution, but this will not be discussed further here.


\subsection{The octonions}
\label{ssOctonions}

The purpose here is to give
an example which does not involve a group and
where the argument given in equation (\ref{eqRiskCalc})
provides a fiducial distribution that
can be used for the determination of the possibility of an
optimal decision rule.
The octonions is here used as an example since it is
one of the more interesting examples of a group-like
structure where the associative law fails.
It has also a natural invariant loss that can be used in the arguments
that follow.
A more familiar example without associativity can be
constructed for the original model
of Fisher for the correlation coefficient,
but we have not been able to identify a natural invariant loss
in that case.

The Cayley--Dickson construction defines a multiplication
$(a,b)(c,d) = (ac-d^* b, d a + b c^*)$ and an involution
$(a,b)^* = (a^*,-b)$ on $\ca\times\ca$ where
$\ca$ is an algebra with an involution.
Starting with the real numbers $\RealN$ this gives
the complex numbers $\CompN$.
Repeated application of the construction gives then
next the quaternions $\QuartN$ and then next the
octonions $\OctN$.
The octonions is hence equal to the $8$-dimensional vector space
$\RealN^8$ equipped with a particular multiplication operation
so that $\OctN$ is an algebra [\citet{Baez02octonions}].

The number $1$ is the unit for multiplication,
and every nonzero element $x$ has a multiplicative inverse $x^{-1}$
with $1 = x x^{-1} = x^{-1} x$.
The usual norm on $\RealN^8$ is\vadjust{\goodbreak} also given by the product and
involution as
$\llVert  x\rrVert^2 = x^* x = x x^*$,
and the identity $\llVert  xy\rrVert  = \llVert  x\rrVert  \llVert
y\rrVert $ holds.
It follows in particular that
$x^{-1} = x^* /\llVert  x\rrVert^2$.
The multiplication is not associative,
but the algebra $\OctN$ is alternative:
the subalgebra generated by any two elements is associative.

Consider next a fiducial model $(U, \chi)$ where
$x = \chi(u,\theta) = \theta u$ is given by the
product in $\OctN$,
and where the conditional law $\pr_U^\theta$ is specified
and does not
depend on $\theta$.
Assume that $\Omega_X = \Omega_U = \Omega_\Theta= \Omega_A = G$
where $G$ is a subset of $\OctN$ that contains $1$,
the product of any two elements, and the inverse of any element.
The particular examples where $G$ is the nonzero octonions
or where $G$ is the octonions with unit norm 
provide examples where $G$ is not a group since the
associative law fails.
This is then a fiducial model for a statistical model
$\{\pr_X^\theta\st\theta\in\Omega_\Theta\}$ where
$(X \st\Theta=\theta) \sim(\theta U \st\Theta= \theta)$.
The corresponding fiducial distribution
is the conditional distribution of
$\Theta^x = x U^{-1}$ given $\Theta= \theta$.
Consider the case where the loss of an action $a$ is given by
$\gamma(\theta, a) = \llVert  \theta- a\rrVert^2 / \llVert  \theta
\rrVert^2$.
This loss is invariant, so the calculation in
equation (\ref{eqRiskCalc}) gives that the risk of
an equivariant decision rule is given by
$\E^\theta\gamma(\Theta^x, \delta(x))$.

Existence of an optimal estimator depends on $\pr_U^\theta$ or
equivalently on $\Theta^x$,
and this will not be discussed further here.
It can, however, be noted that any optimal equivariant decision rule
is determined by $\delta(x) = x \delta(1)$,
and $\delta(1)$ is the minimizer of $\E^\theta\gamma(\Theta^1,
\delta(1))$.
A rule on this form will be equivariant if $\delta(1)$
belongs to the set
$\{a \in G \st(g_1 g_2) a = g_1 (g_2 a)$ $\forall g_1, g_2 \in G\}$.

There are many other examples of binary operations that are not associative.
A~generic family of examples are produced by
a relationship $x = \chi(u, \theta)$ that has the property $(*)$:
it gives a one--one correspondence between the domains of any two of the
variables when the value of the third is fixed.
Corresponding fiducial models based on $\chi$ defines the class of
\textit{simple pivotal models} in accordance with the terminology of
\citet{DawidStone82}, page 1057.
Concrete elementary examples are provided by
$x = u - \theta$, $x = u \theta^{-1}$, and $x = u^\theta$ on
suitable domains.

The property\vspace*{1pt} $(*)$ is conserved by a change of variables by one--one
transformations resulting in
$\phi_x (\tilde{x}) = \chi(\phi_u (\tilde{u}), \phi_\theta
(\tilde{\theta}))$.
For the given elementary examples, there exists
a change of variables so that the result is a relation
$\tilde{x} = \tilde{u} \tilde{\theta}$ given by
a group multiplication.
This is not possible in general.
Simple counter examples arises for the
Fisherian simple pivotal models determined
by the relation $u = F (x \st\theta)$
where $F$ is a suitable cumulative distribution function.
The prototypical example used by Fisher [(\citeyear{Fisher30}), page 534]
when he introduced the fiducial distribution
is given by the sample correlation coefficient from a bivariate
normal distribution.
In this case, a reduction to a group model as for the given elementary
examples is not possible.

In the general case starting from the property $(*)$
there exists, however, a change of variables that results in
a relation given by a quasi-group with a unit: a loop.
The important conclusion of this short discussion is that the
theory of simple pivotal models is linked
naturally to the theory of loops.
The nonzero octonions provides an example of a loop which is not
reducible to a group by a change of variables.

\subsection{Hilbert space}
\label{ssHilbert}

One purpose of this example is to demonstrate existence of
a case where Theorem~\ref{theo1} can be used,
but where an invariant measure does not exist.

Let $\Omega_\Theta= \Omega_A = G$ and $\Omega_X = \Omega_U = G^n$
where $G$ is a complex or real Hilbert space.
The Hilbert space $G$ is a group where the addition of vectors is
the group operation,
and an invariant loss is given by the squared distance between vectors
as $\gamma(\theta, a) = \llVert  \theta- a\rrVert^2$.
A conventional fiducial model $(\chi, U)$ is given by
$x_i = \chi_i (u,\theta) = \theta+ u_i$ for $i=1,\ldots, n$ and a
specification
of a distribution $\pr_U^\theta$ that does not depend on $\theta$.
A maximal invariant is given by
$y = (x_2 - x_1,\ldots, x_n - x_1)$.
The fiducial distribution is given as the
distribution of $\Theta^x = x_1 - U_1$
from the conditional law
$(U \st\Theta=\theta, (U_2-U_1,\ldots, U_n-U_1)=y)$.
The optimal estimator of $\theta$ is given as
$\delta(x) = x_1 - \E(U_1 \st\Theta=\theta, (U_2-U_1,\ldots,U_n-U_1)=y)$.

It will be demonstrated in the next subsection that
it is not necessary to assume independence of $\{U_i\}$
in the previous argument,
and this assumption has indeed not been mentioned above.
More important is the fact that a right Haar prior does not exist
in the case where $G$ is an infinite-dimensional Hilbert space.
An explicit example is given by
$G = l^2 (\NatN) = \{(a_i) \st\llVert  a\rrVert^2 = \sum_{i=1}^\infty\llvert  a_i \rrvert^2 < \infty\}$.

The previous example has an infinite-dimensional parameter space,
and this feature is quite common in applications as exemplified by
nonparametric statistics.
The example does also include data that are infinite dimensional,
and this can be considered to be unrealistic in applications.
There are, however, applications where it is nonetheless convenient
to assume that the observations are also infinite dimensional.
An important source of examples is given by
the statistical signal processing literature [\citet{Trees03ArrayProcessing}].
Explicitly,
it can be convenient to assume that
a signal is observed not only at discretely sampled points,
but for all points.
Similarly, it can be convenient to assume that a complete infinite
sequence of sampled points is observed,
even though only a finite number of samples are actually observed.
In both cases this can lead to a
sample space that is not finite dimensional.
A related and very common convenience is to assume that
observations are given by real numbers,
even though the majority of concrete examples actually only
involves a finite set of observable values due to limited
instrument resolution [\citet{Taraldsen06GUMResolution}].
Explicit consideration of the limit from discretized data
to continuous data gives, incidentally,
a most promising route for the definition of fiducial distributions
more generally than considered in Section~\ref{sOptimal}
as demonstrated recently by \citet{Hannig12fiducial}.

If one observes the real random variables
$X_1,\ldots, X_n$ independently normally distributed
with unknown mean $\theta=(\mu_1,\ldots,\mu_n)$ and variance $1$,\vadjust{\goodbreak}
it is customary to estimate $\mu_i$ by $X_i$.
If the loss is the sum of squares of the errors, this estimator is
admissible for
$n \le2$, but inadmissible for $n > 3$
[\citet{Stein56inadmissible}].
The optimal estimator derived above coincides with the customary estimator.
This exemplifies that the optimal estimator can be inadmissible.
The optimality is only ensured within the class of equivariant estimators.
Equivariance can be a most natural demand,
but this depends on the particular concrete modeling case at hand.
In certain situations [\citet{EfronMorris77steinsparadox}]
it can be natural to give away the equivariance demand in order
to obtain more precise estimates.
In other cases, especially in the context of physics,
the equivariance demand can be closer to the foundation of the subject
matter and will be an absolute demand.



\subsection{Uniform distribution}
\label{ssUniform}

A particular case of the previous Hilbert space example is given
by assuming $G=\RealN$ and where $\pr_U^\theta$ corresponds to a
random sample of size $n$ from the uniform distribution on $(0,1)$.
This gives then a fiducial model for a random sample from the
uniform distribution on $(\theta, \theta+ 1)$.
A fiducial distribution and a corresponding optimal estimator of
$\theta$ follows from the Hilbert space argument.
An alternative and more geometrically tractable argument follows
as explained next
from the use of the sufficient statistic given by the maximum and
minimum observation.

Let $x_i = \theta+ u_i$ where the joint distribution of
$(u_1,u_2)$ conditional on $\Theta=\theta$ is given by
the density $f(u \st\theta) = n (n-1) (u_2 - u_1)^{n-2}$ on
$\{(u_1, u_2) \st0 < u_1 < u_2 < 1\}$.
This is then a fiducial group model for the sufficient statistic
given by the smallest and largest observation from a random sample
from the uniform distribution on $(\theta, \theta+ 1)$.
The model is also a special case of the Hilbert space example with
$n=2$ and where $\{U_i\}$ are conditionally dependent given $\Theta
=\theta$.
Reduction by sufficiency has here simplified the problem,
but the fiducial equation is still over determined so
a further reduction by the maximal invariant
$y = x_2-x_1$ is necessary.
The resulting conditional distribution
$(U_1 \st\Theta=\theta, U_2-U_1=y)$ becomes
the uniform distribution on $(0,1-y)$,
and the fiducial distribution of $\Theta^x$ becomes
the uniform distribution on
$(x_2-1, x_1)$.
This is also a confidence distribution for $\theta$.
The optimal estimator for $\theta$
given the invariant loss $\llvert  \theta- a \rrvert^2$
is $\delta(x) = (x_1 + x_2)/2 - 1/2$.

We choose to add a few comments on this model and estimator
since it has some unusual features.
A first observation is that the Fisher information metric fails to
exist due to nonexistence of the required derivative.
The corresponding distance between two distributions can, however,
still be defined through the length of the parametric
curve $\theta\mapsto\sqrt{f (\cdot\st\theta)}$ in the Hilbert
space of square integrable functions.
This curve is continuous,
but the length from $\theta_1$ to $\theta_2$ is infinite:
it is larger than $2 \sqrt{n} \sqrt{\llvert  \theta_2 - \theta_1
\rrvert }$ for
any integer $n$.

The squared distance $\llvert  \theta_1 - \theta_2 \rrvert^2$ is the squared
distance from the Fisher information metric for any location family
where the density is smooth.
Based on this we consider the invariant loss $\llvert  \theta- a \rrvert^2$
to be a natural choice also in the nonsmooth example considered here.

The optimal estimator $\delta$ found above is unbiased and has hence
minimum variance in the class of unbiased and equivariant estimators.
Nonetheless, according to Lehmann and
Casella [(\citeyear{LehmannCasella98Estimation}), page 87],
there exists no uniformly minimum variance unbiased estimator of
$\theta$.
The statistic $(X_1, X_2)$ is a minimal sufficient statistic,
but it is not complete.
The estimator $\delta$ is, however,
the uniformly minimum variance unbiased estimator
in the larger parametric family which also includes a scale parameter
[\citet{JohnsonKotzBalakrishnanUniv1}, Vol. 2, page 292].
This later reference is also a very good source for further references
and peculiarities regarding the uniform law.\looseness=-1
%

\subsection{Exponential}
\label{sExponential}

The following example is a scale example,
and can be reduced to be a special case of the Hilbert space location
problem by the logarithmic transformation.
A direct solution is equally elementary and is presented to illustrate
the derivation of an optimal estimator.
The explicit formula for the estimator is possibly a novelty.

Assume that $Y_1,\ldots, Y_n$ is a random sample of size $n$ from the
exponential distribution with scale parameter $\beta$.
A fiducial model is given by
$Y_i = \beta V_i$
where the law $\pr_V^\beta$ is as for
a random sample from the standard exponential distribution.
The arithmetic mean $X = \overline{Y}$ is a minimal sufficient statistic.
A corresponding fiducial model is given by
$X = \beta\overline{V} = \beta U$
where $\pr_U^\beta$ is the law of a gamma variable
with scale equal to $1/n$ and shape equal to $n$.
This follows from well-known properties of the gamma distribution.
The model is both simple and conventional,
and the fiducial distribution for an observation $x=\overline{y}$
is hence the conditional distribution of
$\Theta^x = x/U$.
The conclusion is that the fiducial
distribution is the inverse-gamma with
scale $x n$ and shape $n$.

A direct---but more lengthy---calculation of the Bayesian posterior
corresponding to
the right Haar prior $d\beta/ \beta$
gives a posterior that coincides with the fiducial distribution found here.
It is well known more generally that the Bayesian posterior from a
right Haar prior
in a group model coincides with the fiducial.
The calculation demonstrates then that a fiducial model and the
solution of
the fiducial equation gives an alternative and in many cases simpler
route for
the calculation of the Bayesian posterior.
The multivariate normal gives another example where the fiducial
calculation is done in a few lines,
but the corresponding Bayesian calculation is much more cumbersome.

An added advantage of the fiducial calculation is that it shows
directly that the corresponding fiducial distribution is a confidence
distribution.
This is not easily obtained from the Bayesian calculation.
The confidence distribution can alternatively be found
by the likelihood ratio test,
and this has the advantage of giving proof of optimality and
corresponding optimal choices of confidence interval endpoints.
An alternative approach is to also derive optimal intervals
based on Theorem~\ref{theo1} as exemplified by \citet{BERGER}.\vadjust{\goodbreak}

An alternative fiducial calculation can be done without the
reduction to the complete sufficient statistic.
A maximal invariant $\phi$ is given by $\phi(y) = y/\llVert  y\rrVert $.
The conditional law $(V \st\Theta=\theta, \phi(V)=\phi(y))$
will be concentrated on the ray
$Gy = \{\alpha\phi(y) \st\alpha> 0\}$
with a distribution from a density
for $\alpha$ proportional to
$f_V (\alpha\phi(y)) \alpha^{n-1}$.
The assumption of a random sample from the standard exponential
gives a particularly simple $f_V$,
and the fiducial is found explicitly as before.
The alternative calculation has the advantage that it can be
used in the more general case where reduction by sufficiency is not available.


Consider now the problem of estimation of
$\theta=\beta$ with a loss given by
$\gamma(\theta, a) = \llvert  {\ln\theta}- \ln a \rrvert^2$.
This loss is a natural generalization of the squared error loss,
but with the ordinary distance replaced by the distance
$\llvert  {\ln\theta}- \ln a \rrvert $ which is the distance
given by the Fisher information metric in the case of the given scale
model.
In this case,
$G=\Omega_X=\Omega_\Theta=\Omega_A=\RealN^+$
with multiplication as the group operation.
The loss is equivariant, and it follows
that the optimal rule $\delta$ based on the sufficient statistic $X$
is given as the minimizer of
$\rho= \E^\theta\llvert  {\ln\Theta^x} - \ln\delta(x) \rrvert^2$.
This gives that the optimal rule is determined from
$\ln\delta(x) = \E^\theta\ln\Theta^x$.
Evaluation of the corresponding integral
gives an explicit formula for the optimal rule.
It is
%
\begin{equation}
\label{eqScaleLossOpt} \delta(x) = x \exp\bigl(\ln n - \psi(n)\bigr),
\end{equation}
where $\psi$ is the digamma function.
The estimator given by equation (\ref{eqScaleLossOpt})
is possibly known in some contexts,
but we have not found this explicit expression in any of the
textbooks in the list of references or elsewhere.

\subsection{Gamma distribution}
\label{ssGamma}

Assume that $Y_1,\ldots, Y_n$ is a random sample of size $n$ from the
gamma distribution with scale parameter $\beta$ and shape parameter
$\alpha$.
The model parameter is $\theta=(\alpha, \beta)$.
This gives an example as in Section~\ref{ssBernoulli}
where the results related to
Theorem~\ref{theo1} cannot be applied directly.
Fiducial theory can be used to obtain candidates for
good frequentist inference
procedures as indicated next.
Particular results include an exact joint confidence distribution for
$(\alpha, \beta)$,
an exact confidence distribution for $\alpha$,
and a recipe which produces estimators for functions of
$(\alpha, \beta)$ that depends on the choice of a loss.

A fiducial model is given by
$Y_i = \beta F^{-1}(U_i; \alpha)$,
where the law $\pr_U^\theta$ corresponds to a random sample of size $n$
from the uniform distribution on the unit interval $(0,1)$, and $F^{-1}
(u,\alpha)$ is the inverse cumulative distribution function of a gamma
variable with shape $\alpha$ and scale $1$.\vspace*{1pt}

Let $X = (\overline{Y}, \widetilde{Y}/\overline{Y})$ where
$\overline{Y}$ and $\widetilde{Y}$ are the arithmetic and geometric means.
The Bartlett statistic $W=\widetilde{Y}/\overline{Y}$ depends only on
$\alpha$,
and is independent of $\overline{Y}$ as a consequence of the Basu theorem.
A corresponding fiducial model $(\chi, U)$ for $\pr_X^\theta$ is
given by
$\chi_1 (u,\theta) = \beta\overline{F^{-1}(u; \alpha)}$
and
$\chi_2 (u,\theta) = {\widetilde{F^{-1}(u; \alpha)}}/{\overline
{F^{-1}(u; \alpha)}}$.
It\vadjust{\goodbreak} can be noted that $(\chi_2, U)$ gives separately a fiducial model
for $\pr_W^\alpha$.
The corresponding fiducial distribution for $\alpha$ is hence a
confidence distribution.

An alternative fiducial model $(\eta_2, V_2)$ for $\pr_W^\alpha$
is given by inversion of the cumulative distribution function for $W$.
An alternative to $\overline{F^{-1}(u; \alpha)}$ is given by inversion
for a gamma variable with shape $n \alpha$ and scale $1/n$.
The combination gives an alternative fiducial model $(\eta, V)$ for
$\pr_X^\theta$ with the property that $x = \eta(v, \theta)$
defines a one--one correspondence between any two variables
when the third is kept fixed.
The law $\pr_V^\theta$ is the uniform law on the
unit square $(0,1)^2$.
Coordinate transformations can be used to identify
$G = \Omega_\Theta= \Omega_V = \Omega_X$ as sets with a
quasigroup structure with a unit.

Both fiducial models $(\chi, U)$ and $(\eta, V)$
are simple and conventional and determine a fiducial distribution.
For concreteness let $\Theta^x$ be the fiducial corresponding to
$(\eta, V)$.
The quasigroup structure ensures that
$\Theta^x$ gives a joint exact confidence distribution for $(\alpha,
\beta)$.

Consider the problem of estimation of
a function $\tau= h (\alpha, \beta)$ of
the model parameter $\theta= (\alpha, \beta)$.
It can be allowed that $h$ is vector valued,
but assume that each component is positive.
Three examples that are included are then given by
$\tau= \alpha$,
$\tau= \beta$, and
$\tau= (\alpha,\beta)$.
A possible loss in these three cases
is given by squared error loss on a logarithmic scale.
A candidate estimator $\delta$ is then given naturally by
%
\begin{equation}
\label{eqGammaEst1} \delta(x) = \exp \bigl(\E^\theta\ln h\bigl(
\Theta^x\bigr) \bigr).
\end{equation}
This can be evaluated by Monte Carlo simulation from $\pr_V^\theta$
which is the uniform distribution on the unit square $[0,1]^2$.
Another possibility is given by squared distance loss defined
by the Fisher information metric on $\Omega_\Theta$
in the case $h(\theta) = (\alpha, \beta)$.


\section{Discussion}
\label{sDiscussion}

The foundations of Bayesian and frequentist modeling and inference
are well established both in terms of mathematical theory and interpretation.
We do not think that the same can be said about fiducial theory,
but some readers may object to this.
A brief discussion of alternative formulations and naming conventions
found in the literature seems hence to be in place.

Definition~\ref{defFidMod} identifies a fiducial model
with a pair $(U, \zeta)$.
The fiducial model is by definition \textit{conventional} if
$\pr_U^\theta$ does not depend on $\theta$.
In this case we suggest to
denote $U$ as the \textit{Monte Carlo variable} and the
measurable function $\zeta$ as the \textit{fiducial relation}.
The corresponding equation $z = \zeta(u,\theta)$
is the fiducial equation,
but it may also equivalently be denoted
as the fiducial relation.
A fiducial model $(U, \zeta)$ is hence
defined by a Monte Carlo variable $U$ and
a fiducial relation $\zeta$.

If $u$ is a sample from the
Monte Carlo distribution $\pr_U^\theta$,
then $z = \zeta(u, \theta)$ is a sample from the
statistical model as in Definition~\ref{defFidMod}.
The inversion method gives the prototypical
example with $\zeta(u, \theta) = F^{-1}_Z (u \st\theta)$
and $\pr_U^\theta$ equal to the uniform law on the interval $[0,1]$.
This gives the link to the original definition of Fisher,
and also a justification of the choice of the term
\textit{Monte Carlo variable} since this represents
a standard method for simulation from a statistical model on a computer.

The ingredients above given by the pair $(U, \zeta)$
are also the starting point for Dempster--Shafer theory
[\citet{Dempster68bayesgen,Shafer82genbayes}].
\citet{MartinZhangLiu10weakbelief} refer to $U$ as the
\textit{auxiliary variable} and the probability measure
$\mu$ as the \textit{pivotal measure},
where $U \sim\mu$.
The equation $Z = \zeta(U,\Theta)$ is denoted the \textit{a-equation}.
The whole set-up is referred to as an
\textit{inferential model},
and this is identified as something which is not equivalent
to a statistical model.
Except for differences in naming conventions it can be concluded
that an \textit{inferential model} is essentially the same as
a \textit{conventional fiducial model} as summarized in the previous
paragraphs.
The Dempster--Shafer calculus gives an alternative route
for inference based on a fiducial model,
but this is not discussed further here.

The discussion of fiducial theory we have presented is
close to the presentation given by \citet{DawidStone82}.
Dawid and Stone [(\citeyear{DawidStone82}), page 1055] use the term \textit{fiducial model} for
the combination of $\Theta^z = \theta^z (U)$ and
$U \sim\pr_U^\theta$,
and use the term \textit{functional model}
to describe the more general relation
$Z = \zeta(U,\Theta)$.
We chose to avoid the term \textit{functional model} since the term
\textit{functional data analysis} is now the name of a branch of statistics.
\citet{DawidStone82} denote the variable
$U$ as the \textit{error variable},
and uses the symbol $E$ instead.
This corresponds to the naming convention used by \citet{FRASER}.
Fraser [(\citeyear{FRASER}), page 50] uses the terms \textit{structural model} and
\textit{structural equation} in the case
where group theory is an essential ingredient.
\citet{Hannig09fiducial} uses the term
\textit{structural equation} in stead of the term
\textit{fiducial equation} as used by us.
We have avoided the term \textit{structural} here since
there is an active and well-developed different theory
which goes under the label of
\textit{structural equations modeling}.
Our preference for the term \textit{fiducial} as used here
is mainly based on economy of language,
and since this gives the direct link to the original papers of Fisher.

The mathematically inclined reader may claim that
Definition~\ref{defFidMod} is not precise.
This, and the fact that this definition is a novelty
compared with previous writers,
motivate us to state in more detail the assumptions
that are taken as implicitly given from the context in the
statement of Definition~\ref{defFidMod}.
The main difference is that every concept is based
on an underlying abstract conditional probability
space $(\Omega, \ceps, \pr)$ as stated initially
in Section~\ref{sOptimal}.
The fiducial relation $\zeta$ is a measurable function
$\zeta\dvtx  \Omega_U \times\Omega_\Theta\into\Omega_Z$.
This means, as usual,
that $(\zeta\in A) = \{(u,\theta) \st\zeta(u,\theta) \in A\}$
is a measurable set in the product $\sigma$-algebra of
$\Omega_U \times\Omega_\Theta$ whenever $A$ is a measurable
set in $\Omega_Z$.
A consequence is that
$\zeta(U,\Theta)$ is a random element in $\Omega_Z$
defined by the mapping
$\omega\mapsto\zeta(U(\omega),\Theta(\omega))$.
This is measurable since it is assumed
that $\Theta\dvtx  \Omega\into\Omega_\Theta$ and\vadjust{\goodbreak}
$U\dvtx  \Omega\into\Omega_U$ are measurable.
The conditional law $\pr_U^\theta$ of the Monte Carlo variable $U$
is known and does not depend on
$\theta$ in the case of a conventional fiducial model.
If the considerations were limited to the case
where $(\Omega, \ceps, \pr)$ is a probability space,
then this would imply $\pr_U^\theta= \pr_U$.
This fails generally as explained in more detail
by \citet{TaraldsenLindqvist10ImproperPriors}
since $\pr_U^\theta$ is a probability measure,
but $\pr_U$ is unbounded if $\pr$ is unbounded.
The reason for allowing unbounded measures is
the need to include improper priors $\pr_\Theta$.
This has proved fruitful
in related ongoing research by the authors.
It gives in particular
natural conditions that imply equality of
Bayesian posteriors and fiducial distributions.
In specific\vspace*{1pt} modeling cases the spaces
$\Omega_U$, $\Omega_Z$, $\Omega_\Theta$,
the fiducial relation $\zeta$,
and the conditional law $\pr_U^\theta$ are all
explicitly given.
This is as demonstrated by the examples in Section~\ref{sExamples}.
The other ingredients mentioned above are not given
explicitly since they rely on the underlying
unspecified space $\Omega$.
This is as in the ordinary formulation of probability theory
by \citet{KOLMOGOROV} where the whole theory is built
upon the underlying abstract space $\Omega$.
Existence must be proved in each specific modeling case,
but follows trivially in many cases from the construction of
a suitable product space.

Optimal inference for the scale, the location, and the location-scale problems
were investigated using fiducial theory by \citet{Pitman39Location}.
His presentation is most readable and is a good alternative to the
presentations found in standard textbooks.
It can, however, be noted that he concludes that the confidence and
the fiducial theories are essentially the same.
This is in contrast to the views of Neyman and Fisher.
They seemed to agree that in principle the fiducial distribution as
described by Fisher is
not connected to the concept of confidence intervals
as described by Neyman and co-workers.
The content and aims of these two
theories are different.
It seems clear that Fisher never intended to get confidence intervals as
the result of his fiducial arguments.

It is true that the fiducial distributions found in the location-scale
problems, and more general group problem as in
Theorem~\ref{theo1},
are confidence distributions,
but we do consider the concepts to be essentially different in general.
The interpretation of the fiducial distribution,
according to Fisher [(\citeyear{FISHER}), pages 54 and 59]
is identical with the interpretation of the Bayesian posterior:
it represents the state of knowledge regarding the model parameter
as a result of the model assumptions and the observation in the experiment.
It follows then in particular that the fiducial distribution
of a function $\phi(\theta)$ of the model parameter $\theta$
equals the distribution of $\phi(\Theta^x)$ where $\Theta^x$
has the fiducial distribution.
This property does not hold for confidence distributions in general.
In addition,
the fiducial distribution for a simple fiducial model as in
Definition~\ref{defFid1}
is not a confidence distribution in general [\citet{DawidStone82}].

The possibly most famous fiducial distribution is the fiducial
distribution of the difference of means $\mu_1 - \mu_2$ corresponding
to two independent samples from two different normal distributions.
This fiducial distribution gives Fisher's solution to the
Behrens--Fisher problem, but it can be shown by simulation that it is
not a confidence distribution in the sense of having exact coverage
probabilities. A more general class of confidence distributions is
defined by requiring not exact but conservative coverage probabilities.
This is in conformity with the definition of confidence sets in
general. Exactness is often misguidedly taken as a measure of goodness,
but it is not. Power of the associated tests gives one natural measure
of goodness. Examples demonstrate that this more general concept of a
confidence distribution does not coincide with fiducial distributions
in general, but it seems to be an open question whether the
Behrens--Fisher fiducial distribution is a confidence distribution in
this more extended sense. Numerical simulations indicate that the claim
holds
[\citet{Robinson76BehrensFisherAreConservative,Barnard84behrensfisher},
page 269].

The more general problem of obtaining a confidence interval
for the linear combination of several means from
different normal distributions is of considerable practical importance
[\citet{isogum08}].
The ISO recommended solution is in terms of a Welch--Satterthwaite solution,
but a continuation of the arguments given by
\citet{Barnard84behrensfisher} leads to the conclusion that the
fiducial solution is a most competitive alternative solution.

The main virtue of the location-scale models in the context here is
that they illustrate very well the reduction given by a maximal
invariant in cases where a reduction by sufficiency is not possible.
This is also true for the multivariate models treated by
\citet{Fraser79inferenceAndLinear}.
In this case the multivariate normal can be reduced by sufficiency,
but more general models can again be treated by a reduction
through maximal invariants.
It seems that optimal, or good, inference procedures in these
multivariate cases deserves further study guided by fiducial theory.
A recent example of this is given by \citet{LidongHannigIyer08OneWay},
but there are a multitude of different possible examples as indicated by
\citet{Fraser79inferenceAndLinear}.
The suggestion given by Theorem~\ref{theo1} is that not only
confidence intervals,
but also other kinds of inference such as estimation should be considered.

Eaton [(\citeyear{Eaton89groupstat}), pages 89--91] considers the estimation of the
covariance matrix from a multivariate normal sample.
He gives two possible candidates to use as a loss $\gamma(\theta, a)$.
This exemplifies that in the multivariate cases,
and in more complicated group cases, it can
be difficult to decide upon which equivariant loss to use.
It can even be difficult to come up with a good candidate.
In our examples, it has been indicated that the squared distance from
the Fisher information metric is a natural choice.
This will be invariant under mild conditions.
For a statistical model $f (x \st\theta)   \mu(dx)$,
the distance is defined via the length of paths
$t \mapsto x (t) = \sqrt{f (\cdot\st\theta(t))}$ in the Hilbert space
$L_2 (\mu)$.
The nonparametric case given by a parameter space equal to all
densities with respect to $\mu$ gives
the distance $d (f,g) = \cos^{-1} (\int\sqrt{f g}  \,d\mu)$.
The other end of the scale is given by a smooth finite-dimensional
parametric model.
In this case, the previous leads
to the Fisher information metric:
$ds^2 = (1/4) g_{ij} \,d\theta^i \,d\theta^j$ where
$g_{ij} = \E_X^\theta (\partial_i \ln f (X \st\theta) )
(\partial_j \ln f (X \st\theta) )$.
In either case, it gives the model parameter space as
a manifold equipped with a distance
derived intrinsically from the statistical model.


The focus of fiducial theory has initially and currently most often
been on the fiducial distribution by itself and the related possibility
of construction of approximate or exact confidence intervals.
The relation to other kinds of optimal inference such as estimation or
prediction was considered by
Hora and Buehler (\citeyear{HoraBuehler66FiducialInvariant,HoraBuehler67Fiducial}).
The proofs they presented rely on the
existence of an invariant measure,
and it was clear that the fiducial in the case they considered corresponded
to a Bayesian posterior from the right Haar prior.
Since then it has been established in a variety of problems that
the Bayesian algorithm can be used quite generally to obtain good or
optimal frequentist procedures.
The calculation given in equation (\ref{eqRiskCalc}) can be taken as
a strong indication that the fiducial algorithm can be used
similarly to not only obtain confidence intervals,
but also possibly good or
optimal frequentist procedures more generally.
This statement is too general to be provable,
but we consider nonetheless this to be the main content and
message in this paper.
The point of view in this paper does not rely on any particular
interpretation of the fiducial.
It is here simply viewed as a very convenient vehicle for the
derivation of
good, and sometimes optimal as in Theorem~\ref{theo1},
frequentist inference procedures.

\section*{Acknowledgments}

We acknowledge the comments from the reviewers and Associate Editor,
which resulted in an improved paper.



\printaddresses


\begin{thebibliography}{53}

\bibitem[\protect\citeauthoryear{Amari}{1985}]{Amari90diffgeostat}
\begin{bbook}[mr]
\bauthor{\bsnm{Amari},~\bfnm{Shun-ichi}\binits{S.-i.}}
(\byear{1985}).
\btitle{Differential-Geometrical Methods in Statistics}.
\bseries{Lecture Notes in Statistics}
\bvolume{28}.
\bpublisher{Springer}, \blocation{New York}.
\bid{doi={10.1007/978-1-4612-5056-2}, mr={0788689}}
\bptnote{check year}%
\bptok{imsref}%
\end{bbook}
\endbibitem

\bibitem[\protect\citeauthoryear{Anscombe}{1948}]{Anscombe48fiducial}
\begin{barticle}[mr]
\bauthor{\bsnm{Anscombe},~\bfnm{F.~J.}\binits{F.~J.}}
(\byear{1948}).
\btitle{The validity of comparative experiments}.
\bjournal{J. Roy. Statist. Soc. Ser. A}
\bvolume{111}
\bpages{181--200; discussion, 200--211}.
\bid{issn={0035-9238}, mr={0030181}}
\bptnote{check related}%
\bptok{imsref}%
\end{barticle}
\endbibitem

\bibitem[\protect\citeauthoryear{Atkinson and
Mitchell}{1981}]{AtkinsonMitchell81statdistance}
\begin{barticle}[mr]
\bauthor{\bsnm{Atkinson},~\bfnm{Colin}\binits{C.}} \AND
\bauthor{\bsnm{Mitchell},~\bfnm{Ann F.~S.}\binits{A.~F.~S.}}
(\byear{1981}).
\btitle{Rao's distance measure}.
\bjournal{Sankhy\=a Ser. A}
\bvolume{43}
\bpages{345--365}.
\bid{issn={0581-572X}, mr={0665876}}
\bptok{imsref}%
\end{barticle}
\endbibitem

\bibitem[\protect\citeauthoryear{Baez}{2002}]{Baez02octonions}
\begin{barticle}[mr]
\bauthor{\bsnm{Baez},~\bfnm{John~C.}\binits{J.~C.}}
(\byear{2002}).
\btitle{The octonions}.
\bjournal{Bull. Amer. Math. Soc. (N.S.)}
\bvolume{39}
\bpages{145--205}.
\bid{doi={10.1090/S0273-0979-01-00934-X}, issn={0273-0979}, mr={1886087}}
\bptok{imsref}%
\end{barticle}
\endbibitem

\bibitem[\protect\citeauthoryear{Barnard}{1984}]{Barnard84behrensfisher}
\begin{barticle}[mr]
\bauthor{\bsnm{Barnard},~\bfnm{G.~A.}\binits{G.~A.}}
(\byear{1984}).
\btitle{Comparing the means of two independent samples}.
\bjournal{J. R. Stat. Soc. Ser. C. Appl. Stat.}
\bvolume{33}
\bpages{266--271}.
\bid{doi={10.2307/2347702}, issn={0035-9254}, mr={0782073}}
\bptok{imsref}%
\end{barticle}
\endbibitem

\bibitem[\protect\citeauthoryear{Berger}{1985}]{BERGER}
\begin{bbook}[mr]
\bauthor{\bsnm{Berger},~\bfnm{James~O.}\binits{J.~O.}}
(\byear{1985}).
\btitle{Statistical Decision Theory and {B}ayesian Analysis},
\bedition{2nd} ed.
\bpublisher{Springer}, \blocation{New York}.
\bid{mr={0804611}}
\bptok{imsref}%
\end{bbook}
\endbibitem

\bibitem[\protect\citeauthoryear{Blank}{1956}]{Blank56binomialoptimal}
\begin{barticle}[mr]
\bauthor{\bsnm{Blank},~\bfnm{A.~A.}\binits{A.~A.}}
(\byear{1956}).
\btitle{Existence and uniqueness of a uniformly most powerful randomized
unbiased test for the binomial}.
\bjournal{Biometrika}
\bvolume{43}
\bpages{465--466}.
\bid{issn={0006-3444}, mr={0081601}}
\bptok{imsref}%
\end{barticle}
\endbibitem

\bibitem[\protect\citeauthoryear{Bunke}{1975}]{Bunke75functional}
\begin{barticle}[mr]
\bauthor{\bsnm{Bunke},~\bfnm{H.}\binits{H.}}
(\byear{1975}).
\btitle{Statistical inference: Fiducial and structural vs. likelihood}.
\bjournal{Math. Operationsforsch. Statist.}
\bvolume{6}
\bpages{667--676}.
\bid{mr={0394946}}
\bptok{imsref}%
\end{barticle}
\endbibitem

\bibitem[\protect\citeauthoryear{Castaing and Valadier}{1977}]{CASTAING}
\begin{bbook}[mr]
\bauthor{\bsnm{Castaing},~\bfnm{C.}\binits{C.}} \AND
\bauthor{\bsnm{Valadier},~\bfnm{M.}\binits{M.}}
(\byear{1977}).
\btitle{Convex Analysis and Measurable Multifunctions}.
\bseries{Lecture Notes in Math.}
\bvolume{580}.
\bpublisher{Springer}, \blocation{Berlin}.
\bid{mr={0467310}}
\bptok{imsref}%
\end{bbook}
\endbibitem

\bibitem[\protect\citeauthoryear{Dawid and Stone}{1982}]{DawidStone82}
\begin{barticle}[mr]
\bauthor{\bsnm{Dawid},~\bfnm{A.~P.}\binits{A.~P.}} \AND
\bauthor{\bsnm{Stone},~\bfnm{M.}\binits{M.}}
(\byear{1982}).
\btitle{The functional-model basis of fiducial inference (with discussion)}.
\bjournal{Ann. Statist.}
\bvolume{10}
\bpages{1054--1074}.
\bid{issn={0090-5364}, mr={0673643}}
\bptok{imsref}%
\end{barticle}
\endbibitem

\bibitem[\protect\citeauthoryear{Dempster}{1968}]{Dempster68bayesgen}
\begin{barticle}[mr]
\bauthor{\bsnm{Dempster},~\bfnm{A.~P.}\binits{A.~P.}}
(\byear{1968}).
\btitle{A generalization of {B}ayesian inference (with discussion)}.
\bjournal{J. R. Stat. Soc. Ser. B Stat. Methodol.}
\bvolume{30}
\bpages{205--247}.
\bid{issn={0035-9246}, mr={0238428}}
\bptok{imsref}%
\end{barticle}
\endbibitem

\bibitem[\protect\citeauthoryear{E, Hannig and
Iyer}{2008}]{LidongHannigIyer08OneWay}
\begin{barticle}[mr]
\bauthor{\bsnm{E},~\bfnm{Lidong}\binits{L.}},
\bauthor{\bsnm{Hannig},~\bfnm{Jan}\binits{J.}} \AND
\bauthor{\bsnm{Iyer},~\bfnm{Hari}\binits{H.}}
(\byear{2008}).
\btitle{Fiducial intervals for variance components in an unbalanced
two-component normal mixed linear model}.
\bjournal{J. Amer. Statist. Assoc.}
\bvolume{103}
\bpages{854--865}.
\bid{doi={10.1198/016214508000000229}, issn={0162-1459}, mr={2524335}}
\bptok{imsref}%
\end{barticle}
\endbibitem

\bibitem[\protect\citeauthoryear{Eaton}{1989}]{Eaton89groupstat}
\begin{bbook}[auto:STB|2013/01/29|08:09:18]
\bauthor{\bsnm{Eaton},~\bfnm{M.~L.}\binits{M.~L.}}
(\byear{1989}).
\btitle{Group Invariance Applications in Statistics}
\bvolume{1}.
\bpublisher{IMS}, \blocation{Hayward, CA}.
\bptok{imsref}%
\end{bbook}
\endbibitem

\bibitem[\protect\citeauthoryear{Efron}{1998}]{Efron98}
\begin{barticle}[mr]
\bauthor{\bsnm{Efron},~\bfnm{Bradley}\binits{B.}}
(\byear{1998}).
\btitle{R. {A}. {F}isher in the 21st century (with discussion)}.
\bjournal{Statist. Sci.}
\bvolume{13}
\bpages{95--122}.
\bid{doi={10.1214/ss/1028905930}, issn={0883-4237}, mr={1647499}}
\bptok{imsref}%
\end{barticle}
\endbibitem


\bibitem[\protect\citeauthoryear{Efron}{2006}]{Efron06fiducial}
\begin{barticle}[mr]
\bauthor{\bsnm{Efron},~\bfnm{Bradley}\binits{B.}}
(\byear{2006}).
\btitle{Minimum volume confidence regions for a multivariate normal mean
vector}.
\bjournal{J.~R. Stat. Soc. Ser. B Stat. Methodol.}
\bvolume{68}
\bpages{655--670}.
\bid{doi={10.1111/j.1467-9868.2006.00560.x}, issn={1369-7412}, mr={2301013}}
\bptok{imsref}%
\end{barticle}
\endbibitem

\bibitem[\protect\citeauthoryear{Efron and
Morris}{1977}]{EfronMorris77steinsparadox}
\begin{barticle}[auto:STB|2013/01/29|08:09:18]
\bauthor{\bsnm{Efron},~\bfnm{B.}\binits{B.}} \AND
\bauthor{\bsnm{Morris},~\bfnm{C.}\binits{C.}}
(\byear{1977}).
\btitle{Stein's paradox in statistics}.
\bjournal{Scientific American}
\bvolume{236}
\bpages{119--127}.
\bptok{imsref}%
\end{barticle}
\endbibitem

\bibitem[\protect\citeauthoryear{Fisher}{1930}]{Fisher30}
\begin{barticle}[auto:STB|2013/01/29|08:09:18]
\bauthor{\bsnm{Fisher},~\bfnm{R.~A.}\binits{R.~A.}}
(\byear{1930}).
\btitle{Inverse probability}.
\bjournal{Math. Proc. Cambridge Philos. Soc.}
\bvolume{26}
\bpages{528--535}.
\bptok{imsref}%
\end{barticle}
\endbibitem

\bibitem[\protect\citeauthoryear{Fisher}{1935}]{Fisher35Fiducial}
\begin{barticle}[auto:STB|2013/01/29|08:09:18]
\bauthor{\bsnm{Fisher},~\bfnm{R.~A.}\binits{R.~A.}}
(\byear{1935}).
\btitle{The fiducial argument in statistical inference}.
\bjournal{Annals of Eugenics}
\bvolume{6}
\bpages{391--398}.
\bptok{imsref}%
\end{barticle}
\endbibitem

\bibitem[\protect\citeauthoryear{Fisher}{1973}]{FISHER}
\begin{bbook}[mr]
\bauthor{\bsnm{Fisher},~\bfnm{Ronald~A.}\binits{R.~A.}}
(\byear{1973}).
\btitle{Statistical Methods and Scientific Inference}.
\bpublisher{Hafner}, \blocation{New York}.
\bptok{imsref}%
\end{bbook}
\endbibitem

\bibitem[\protect\citeauthoryear{Fraser}{1961a}]{Fraser61Fiducial}
\begin{barticle}[mr]
\bauthor{\bsnm{Fraser},~\bfnm{D.~A.~S.}\binits{D.~A.~S.}}
(\byear{1961}a).
\btitle{On fiducial inference}.
\bjournal{Ann. Math. Statist.}
\bvolume{32}
\bpages{661--676}.
\bid{issn={0003-4851}, mr={0130755}}
\bptok{imsref}%
\end{barticle}
\endbibitem

\bibitem[\protect\citeauthoryear{Fraser}{1961b}]{Fraser61FiducialInvariance}
\begin{barticle}[mr]
\bauthor{\bsnm{Fraser},~\bfnm{D.~A.~S.}\binits{D.~A.~S.}}
(\byear{1961}b).
\btitle{The fiducial method and invariance}.
\bjournal{Biometrika}
\bvolume{48}
\bpages{261--280}.
\bid{issn={0006-3444}, mr={0133910}}
\bptok{imsref}%
\end{barticle}
\endbibitem

\bibitem[\protect\citeauthoryear{Fraser}{1962}]{Fraser62fiducial}
\begin{barticle}[auto:STB|2013/01/29|08:09:18]
\bauthor{\bsnm{Fraser},~\bfnm{D.~A.~S.}\binits{D.~A.~S.}}
(\byear{1962}).
\btitle{On the consistency of the fiducial method}.
\bjournal{J. R. Stat. Soc. Ser. B Stat. Methodol.}
\bvolume{24}
\bpages{425--434}.
\bptok{imsref}%
\end{barticle}
\endbibitem

\bibitem[\protect\citeauthoryear{Fraser}{1963}]{Fraser64fiducial}
\begin{barticle}[mr]
\bauthor{\bsnm{Fraser},~\bfnm{D.~A.~S.}\binits{D.~A.~S.}}
(\byear{1963}).
\btitle{On the definition of fiducial probability}.
\bjournal{Bull. Inst. Internat. Statist.}
\bvolume{40}
\bpages{842--856}.
\bid{mr={0172439}}
\bptnote{check year}%
\bptok{imsref}%
\end{barticle}
\endbibitem

\bibitem[\protect\citeauthoryear{Fraser}{1968}]{FRASER}
\begin{bbook}[mr]
\bauthor{\bsnm{Fraser},~\bfnm{D.~A.~S.}\binits{D.~A.~S.}}
(\byear{1968}).
\btitle{The Structure of Inference}.
\bpublisher{Wiley}, \blocation{New York}.
\bid{mr={0235643}}
\bptok{imsref}%
\end{bbook}
\endbibitem

\bibitem[\protect\citeauthoryear{Fraser}{1979}]{Fraser79inferenceAndLinear}
\begin{bbook}[mr]
\bauthor{\bsnm{Fraser},~\bfnm{Donald Alexander~Stuart}\binits{D.~A.~S.}}
(\byear{1979}).
\btitle{Inference and Linear Models: Advanced Book Program}.
\bpublisher{McGraw-Hill}, \blocation{New York}.
\bid{mr={0535612}}
\bptok{imsref}%
\end{bbook}
\endbibitem

\bibitem[\protect\citeauthoryear{Fraser
et~al.}{2010}]{FraserReidMarrasYi10fiducial}
\begin{barticle}[mr]
\bauthor{\bsnm{Fraser},~\bfnm{D.~A.~S.}\binits{D.~A.~S.}},
\bauthor{\bsnm{Reid},~\bfnm{N.}\binits{N.}},
\bauthor{\bsnm{Marras},~\bfnm{E.}\binits{E.}} \AND
\bauthor{\bsnm{Yi},~\bfnm{G.~Y.}\binits{G.~Y.}}
(\byear{2010}).
\btitle{Default priors for {B}ayesian and frequentist inference}.
\bjournal{J. R. Stat. Soc. Ser. B Stat. Methodol.}
\bvolume{72}
\bpages{631--654}.
\bid{doi={10.1111/j.1467-9868.2010.00750.x}, issn={1369-7412}, mr={2758239}}
\bptok{imsref}%
\end{barticle}
\endbibitem

\bibitem[\protect\citeauthoryear{Ghosh, Reid and
Fraser}{2010}]{GhoshReidFraser10CondInference}
\begin{barticle}[mr]
\bauthor{\bsnm{Ghosh},~\bfnm{M.}\binits{M.}},
\bauthor{\bsnm{Reid},~\bfnm{N.}\binits{N.}} \AND
\bauthor{\bsnm{Fraser},~\bfnm{D.~A.~S.}\binits{D.~A.~S.}}
(\byear{2010}).
\btitle{Ancillary statistics: A review}.
\bjournal{Statist. Sinica}
\bvolume{20}
\bpages{1309--1332}.
\bid{issn={1017-0405}, mr={2777327}}
\bptok{imsref}%
\end{barticle}
\endbibitem

\bibitem[\protect\citeauthoryear{Hannig}{2009}]{Hannig09fiducial}
\begin{barticle}[mr]
\bauthor{\bsnm{Hannig},~\bfnm{Jan}\binits{J.}}
(\byear{2009}).
\btitle{On generalized fiducial inference}.
\bjournal{Statist. Sinica}
\bvolume{19}
\bpages{491--544}.
\bid{issn={1017-0405}, mr={2514173}}
\bptok{imsref}%
\end{barticle}
\endbibitem

\bibitem[\protect\citeauthoryear{Hannig}{2013}]{Hannig12fiducial}
\begin{bmisc}[auto:STB|2013/01/29|08:09:18]
\bauthor{\bsnm{Hannig},~\bfnm{J.}\binits{J.}}
(\byear{2013}).
\bhowpublished{Generalized fiducial inference via discretization.
\textit{Statist. Sinica}. To appear}.
\bptok{imsref}%
\end{bmisc}
\endbibitem

\bibitem[\protect\citeauthoryear{Hora and
Buehler}{1966}]{HoraBuehler66FiducialInvariant}
\begin{barticle}[mr]
\bauthor{\bsnm{Hora},~\bfnm{R.~B.}\binits{R.~B.}} \AND
\bauthor{\bsnm{Buehler},~\bfnm{R.~J.}\binits{R.~J.}}
(\byear{1966}).
\btitle{Fiducial theory and invariant estimation}.
\bjournal{Ann. Math. Statist.}
\bvolume{37}
\bpages{643--656}.
\bid{issn={0003-4851}, mr={0199938}}
\bptok{imsref}%
\end{barticle}
\endbibitem

\bibitem[\protect\citeauthoryear{Hora and
Buehler}{1967}]{HoraBuehler67Fiducial}
\begin{barticle}[mr]
\bauthor{\bsnm{Hora},~\bfnm{R.~B.}\binits{R.~B.}} \AND
\bauthor{\bsnm{Buehler},~\bfnm{R.~J.}\binits{R.~J.}}
(\byear{1967}).
\btitle{Fiducial theory and invariant prediction}.
\bjournal{Ann. Math. Statist.}
\bvolume{38}
\bpages{795--801}.
\bid{issn={0003-4851}, mr={0211520}}
\bptok{imsref}%
\end{barticle}
\endbibitem

\bibitem[\protect\citeauthoryear{ISO/IEC}{2008}]{isogum08}
\begin{bmisc}[auto:STB|2013/01/29|08:09:18]
\borganization{ISO/IEC}
(\byear{2008}).
\bhowpublished{Guide 98-3: Uncertainty of measurement---Part 3: Guide to the
expression of uncertainty in measurement (GUM:1995). Technical report,
Geneve}.
\bptok{imsref}%
\end{bmisc}
\endbibitem

\bibitem[\protect\citeauthoryear{Johnson, Kotz and
Balakrishnan}{1994}]{JohnsonKotzBalakrishnanUniv1}
\begin{bbook}[mr]
\bauthor{\bsnm{Johnson},~\bfnm{Norman~L.}\binits{N.~L.}},
\bauthor{\bsnm{Kotz},~\bfnm{Samuel}\binits{S.}} \AND
\bauthor{\bsnm{Balakrishnan},~\bfnm{N.}\binits{N.}}
(\byear{1994}).
\btitle{Continuous Univariate Distributions. {V}ol. 1},
\bedition{2nd} ed.
\bpublisher{Wiley}, \blocation{New York}.
\bid{mr={1299979}}
\bptok{imsref}%
\end{bbook}
\endbibitem

\bibitem[\protect\citeauthoryear{Kolmogorov}{1956}]{KOLMOGOROV}
\begin{bbook}[mr]
\bauthor{\bsnm{Kolmogorov},~\bfnm{A.~N.}\binits{A.~N.}}
(\byear{1956}).
\btitle{Foundations of the Theory of Probability}.
\bpublisher{Chelsea}, \blocation{New York}.
\bid{mr={0079843}}
\bptok{imsref}%
\end{bbook}
\endbibitem

\bibitem[\protect\citeauthoryear{Lehmann and
Casella}{1998}]{LehmannCasella98Estimation}
\begin{bbook}[mr]
\bauthor{\bsnm{Lehmann},~\bfnm{E.~L.}\binits{E.~L.}} \AND
\bauthor{\bsnm{Casella},~\bfnm{George}\binits{G.}}
(\byear{1998}).
\btitle{Theory of Point Estimation},
\bedition{2nd} ed.
\bpublisher{Springer}, \blocation{New York}.
\bid{mr={1639875}}
\bptok{imsref}%
\end{bbook}
\endbibitem

\bibitem[\protect\citeauthoryear{Lehmann and
Romano}{2005}]{LehmannRomano05testing}
\begin{bbook}[mr]
\bauthor{\bsnm{Lehmann},~\bfnm{E.~L.}\binits{E.~L.}} \AND
\bauthor{\bsnm{Romano},~\bfnm{Joseph~P.}\binits{J.~P.}}
(\byear{2005}).
\btitle{Testing Statistical Hypotheses},
\bedition{3rd} ed.
\bpublisher{Springer}, \blocation{New York}.
\bid{mr={2135927}}
\bptok{imsref}%
\end{bbook}
\endbibitem

\bibitem[\protect\citeauthoryear{Martin, Zhang and
Liu}{2010}]{MartinZhangLiu10weakbelief}
\begin{barticle}[mr]
\bauthor{\bsnm{Martin},~\bfnm{Ryan}\binits{R.}},
\bauthor{\bsnm{Zhang},~\bfnm{Jianchun}\binits{J.}} \AND
\bauthor{\bsnm{Liu},~\bfnm{Chuanhai}\binits{C.}}
(\byear{2010}).
\btitle{Dempster--{S}hafer theory and statistical inference with weak beliefs}.
\bjournal{Statist. Sci.}
\bvolume{25}
\bpages{72--87}.
\bid{doi={10.1214/10-STS322}, issn={0883-4237}, mr={2741815}}
\bptok{imsref}%
\end{barticle}
\endbibitem

\bibitem[\protect\citeauthoryear{Pitman}{1939}]{Pitman39Location}
\begin{barticle}[auto:STB|2013/01/29|08:09:18]
\bauthor{\bsnm{Pitman},~\bfnm{E.~J.~G.}\binits{E.~J.~G.}}
(\byear{1939}).
\btitle{The estimation of the location and scale parameters of a continuous
population of any given form}.
\bjournal{Biometrika}
\bvolume{30}
\bpages{391--421}.
\bptok{imsref}%
\end{barticle}
\endbibitem

\bibitem[\protect\citeauthoryear{Radhakrishna~Rao}{1945}]{Rao45}
\begin{barticle}[mr]
\bauthor{\bsnm{Radhakrishna~Rao},~\bfnm{C.}\binits{C.}}
(\byear{1945}).
\btitle{Information and the accuracy attainable in the estimation of
statistical parameters}.
\bjournal{Bull. Calcutta Math. Soc.}
\bvolume{37}
\bpages{81--91}.
\bid{issn={0008-0659}, mr={0015748}}
\bptok{imsref}%
\end{barticle}
\endbibitem

\bibitem[\protect\citeauthoryear{Robinson}{1976}]{Robinson76BehrensFisherAreCo%
nservative}
\begin{barticle}[mr]
\bauthor{\bsnm{Robinson},~\bfnm{G.~K.}\binits{G.~K.}}
(\byear{1976}).
\btitle{Properties of {S}tudent's {$t$} and of the {B}ehrens--{F}isher solution
to the two means problem}.
\bjournal{Ann. Statist.}
\bvolume{4}
\bpages{963--971}.
\bid{issn={0090-5364}, mr={0415868}}
\bptok{imsref}%
\end{barticle}
\endbibitem

\bibitem[\protect\citeauthoryear{Schervish}{1995}]{SCHERVISH}
\begin{bbook}[mr]
\bauthor{\bsnm{Schervish},~\bfnm{Mark~J.}\binits{M.~J.}}
(\byear{1995}).
\btitle{Theory of Statistics}.
\bpublisher{Springer}, \blocation{New York}.
\bid{doi={10.1007/978-1-4612-4250-5}, mr={1354146}}
\bptok{imsref}%
\end{bbook}
\endbibitem

\bibitem[\protect\citeauthoryear{Shafer}{1982}]{Shafer82genbayes}
\begin{barticle}[mr]
\bauthor{\bsnm{Shafer},~\bfnm{Glenn}\binits{G.}}
(\byear{1982}).
\btitle{Belief functions and parametric models (with discussion)}.
\bjournal{J. R. Stat. Soc. Ser. B Stat. Methodol.}
\bvolume{44}
\bpages{322--352}.
\bid{issn={0035-9246}, mr={0693232}}
\bptok{imsref}%
\end{barticle}
\endbibitem

\bibitem[\protect\citeauthoryear{Smith}{2007}]{Smith06loop}
\begin{bbook}[mr]
\bauthor{\bsnm{Smith},~\bfnm{Jonathan D.~H.}\binits{J.~D.~H.}}
(\byear{2007}).
\btitle{An Introduction to Quasigroups and Their Representations}.
\bpublisher{Chapman \& Hall/CRC}, \blocation{Boca Raton, FL}.
\bid{mr={2268350}}
\bptnote{check year}%
\bptok{imsref}%
\end{bbook}
\endbibitem

\bibitem[\protect\citeauthoryear{Stein}{1956}]{Stein56inadmissible}
\begin{binproceedings}[mr]
\bauthor{\bsnm{Stein},~\bfnm{Charles}\binits{C.}}
(\byear{1956}).
\btitle{Inadmissibility of the usual estimator for the mean of a multivariate
normal distribution}.
In \bbooktitle{Proceedings of the {T}hird {B}erkeley {S}ymposium on
{M}athematical {S}tatistics and {P}robability, 1954--1955, Vol. {I}}
\bpages{197--206}.
\bpublisher{Univ. California Press}, \blocation{Berkeley and Los Angeles}.
\bid{mr={0084922}}
\bptok{imsref}%
\end{binproceedings}
\endbibitem

\bibitem[\protect\citeauthoryear{Stevens}{1950}]{Stevens50binomial}
\begin{barticle}[mr]
\bauthor{\bsnm{Stevens},~\bfnm{W.~L.}\binits{W.~L.}}
(\byear{1950}).
\btitle{Fiducial limits of the parameter of a discontinuous distribution}.
\bjournal{Biometrika}
\bvolume{37}
\bpages{117--129}.
\bid{issn={0006-3444}, mr={0035955}}
\bptok{imsref}%
\end{barticle}
\endbibitem

\bibitem[\protect\citeauthoryear{Stevens}{1957}]{Stevens57binomial}
\begin{barticle}[mr]
\bauthor{\bsnm{Stevens},~\bfnm{W.~L.}\binits{W.~L.}}
(\byear{1957}).
\btitle{Shorter intervals for the parameter of the binomial and {P}oisson
distributions}.
\bjournal{Biometrika}
\bvolume{44}
\bpages{436--440}.
\bid{issn={0006-3444}, mr={0090188}}
\bptok{imsref}%
\end{barticle}
\endbibitem

\bibitem[\protect\citeauthoryear{Stuart, Ord and
Arnold}{1999}]{StuartOrdArnold99kendal2a}
\begin{bbook}[auto:STB|2013/01/29|08:09:18]
\bauthor{\bsnm{Stuart},~\bfnm{A.}\binits{A.}},
\bauthor{\bsnm{Ord},~\bfnm{K.}\binits{K.}} \AND
\bauthor{\bsnm{Arnold},~\bfnm{S.}\binits{S.}}
(\byear{1999}).
\btitle{Kendall's Advanced Theory of Statistics, Classical Inference and the
Linear Model (2007 Reprint)},
\bedition{6th} ed.
\bvolume{2A}.
\bpublisher{Wiley}, \blocation{New York}.
\bptok{imsref}%
\end{bbook}
\endbibitem

\bibitem[\protect\citeauthoryear{Taraldsen}{2006}]{Taraldsen06GUMResolution}
\begin{barticle}[auto:STB|2013/01/29|08:09:18]
\bauthor{\bsnm{Taraldsen},~\bfnm{G.}\binits{G.}}
(\byear{2006}).
\btitle{Instrument resolution and measurement accuracy}.
\bjournal{Metrologia}
\bvolume{43}
\bpages{539--544}.
\bptok{imsref}%
\end{barticle}
\endbibitem

\bibitem[\protect\citeauthoryear{Taraldsen and
Lindqvist}{2010}]{TaraldsenLindqvist10ImproperPriors}
\begin{barticle}[mr]
\bauthor{\bsnm{Taraldsen},~\bfnm{Gunnar}\binits{G.}} \AND
\bauthor{\bsnm{Lindqvist},~\bfnm{Bo~Henry}\binits{B.~H.}}
(\byear{2010}).
\btitle{Improper priors are not improper}.
\bjournal{Amer. Statist.}
\bvolume{64}
\bpages{154--158}.
\bid{doi={10.1198/tast.2010.09116}, issn={0003-1305}, mr={2757006}}
\bptok{imsref}%
\end{barticle}
\endbibitem

\bibitem[\protect\citeauthoryear{Van~Trees}{2003}]{Trees03ArrayProcessing}
\begin{bbook}[auto:STB|2013/01/29|08:09:18]
\bauthor{\bsnm{Van~Trees},~\bfnm{H.~L.}\binits{H.~L.}}
(\byear{2003}).
\btitle{Detection, Estimation, and Modulation Theory (Volumes:~\mbox{I--IV})}.
\bpublisher{Wiley}, \blocation{New York}.
\bptok{imsref}%
\end{bbook}
\endbibitem

\bibitem[\protect\citeauthoryear{Wang, Hannig and
Iyer}{2012}]{WangHannigIyer12gumfiducial}
\begin{barticle}[auto:STB|2013/01/29|08:09:18]
\bauthor{\bsnm{Wang},~\bfnm{C.~M.}\binits{C.~M.}},
\bauthor{\bsnm{Hannig},~\bfnm{J.}\binits{J.}} \AND
\bauthor{\bsnm{Iyer},~\bfnm{H.~K.}\binits{H.~K.}}
(\byear{2012}).
\btitle{Pivotal methods in the propagation of distributions}.
\bjournal{Metrologia}
\bvolume{49}
\bpages{382--389}.
\bptok{imsref}%
\end{barticle}
\endbibitem

\bibitem[\protect\citeauthoryear{Wijsman}{1990}]{Wijsman90statgroups}
\begin{bbook}[mr]
\bauthor{\bsnm{Wijsman},~\bfnm{Robert~A.}\binits{R.~A.}}
(\byear{1990}).
\btitle{Invariant Measures on Groups and Their Use in Statistics}.
\bseries{Institute of Mathematical Statistics Lecture Notes---Monograph Series}
\bvolume{14}.
\bpublisher{IMS}, \blocation{Hayward, CA}.
\bid{mr={1218397}}
\bptok{imsref}%
\end{bbook}
\endbibitem

\end{thebibliography}
\end{document}